\documentclass[twocolumn]{autart}    

\usepackage{graphicx}          

\usepackage{amsmath}
\usepackage{amssymb}

\newcommand*\diff{\mathop{}\!\mathrm{d}}%
\DeclareMathOperator*{\esssup}{ess\,sup}

\usepackage{graphics} 
\usepackage{epsfig} 
\usepackage{epstopdf}
\usepackage{subfigure}
\usepackage{color}

\graphicspath{{./Figures/}}

\allowdisplaybreaks[4]

\begin{document}

\begin{frontmatter}

\title{Boundary feedback stabilization of a reaction-diffusion equation with Robin boundary conditions and state-delay\thanksref{footnoteinfo}} 

\thanks[footnoteinfo]{This work was supported by a research grant from Science Foundation Ireland (SFI) under grant number 16/RC/3872 and is co-funded under the European Regional Development Fund and by I-Form industry partners.\\ 
Corresponding author H.~Lhachemi.}

\author[UCD]{Hugo Lhachemi}\ead{hugo.lhachemi@ucd.ie},
\author[UCD]{Robert Shorten}\ead{robert.shorten@ucd.ie},               

\address[UCD]{School of Electrical and Electronic Engineering, University College Dublin, Dublin, Ireland}  


\begin{abstract}                          
This paper discusses the boundary feedback stabilization of a reaction-diffusion equation with Robin boundary conditions and in the presence of a time-varying state-delay. The proposed control design strategy is based on a finite-dimensional truncated model obtained via a spectral decomposition. By an adequate selection of the number of modes of the original infinite-dimensional system, we show that the design performed on the finite-dimensional truncated model achieves the exponential stabilization of the original infinite-dimensional system. In the presence of distributed disturbances, we show that the closed-loop system is exponentially input-to-state stable with fading memory. 
\end{abstract}

\end{frontmatter}

\section{Introduction}
Boundary stabilization of partial differential equations (PDEs) in the presence of delays is an active research topic. A first research direction deals with the feedback stabilization of PDEs by means of delayed boundary control. For example, the cases of the heat~\cite{nicaise2009stability} and wave~\cite{nicaise2008stabilization,nicaise2007stabilization,nicaise2009stability} equations were studied via Lyapunov methods for slowly time-varying delays. A backstepping approach was reported in~\cite{krstic2009control} for the boundary feedback stabilization of an unstable reaction-diffusion equation under large constant input delays. Inspired by the early work~\cite{russell1978controllability} and the developments reported in~\cite{coron2004global,coron2006global}, such a problem was also investigated in~\cite{prieur2018feedback} by designing a predictor feedback on a finite-dimensional truncated model capturing the unstable modes of the system. Then, a Lyapunov-based argument was employed to ensure that the control law achieves the stabilization of the full infinite-dimensional system. The same approach was applied in~\cite{guzman2019stabilization} for the boundary stabilization of a linear Kuramoto-Sivashinsky equation. Such an approach was generalized to a class of diagonal infinite-dimensional systems in~\cite{lhachemi2019feedback,lhachemi2019control} for constant input delays and then in~\cite{lhachemi2019lmi} for fast time-varying input delays. A second research direction deals with the boundary feedback stabilization of PDEs in the presence of a state-delay. Such state-delays can be used to model either/both locality or/and inertia of certain physical phenomenon such as heat or mass transfers~\cite{polyanin2014nonlinear}. Motivated by the success of Linear Matrix Inequalities (LMI)-based approaches for the study of delayed finite-dimensional systems~\cite{fridman2014tutorial}, LMI conditions were investigated in~\cite{fridman2009exponential,solomon2015stability} for the stability analysis of PDEs in the presence of a state-delay. For the boundary control design of state-delayed PDEs, backstepping-based methods were reported in~\cite{hashimoto2016stabilization,kang2017boundary,kang2017boundaryIFAC,kang2018boundary}.

This paper is concerned with the boundary feedback stabilization of a reaction-diffusion equation with Robin boundary conditions in the presence of a time-varying state-delay. A similar setting was investigated in~\cite{hashimoto2016stabilization} via a backstepping approach in the case of Dirichlet-Dirichlet boundary conditions and for a constant state-delay. This problem was also investigated by means of backstepping control design in~\cite{kang2017boundaryIFAC} for Neumann-Dirichlet boundary conditions with Dirichlet actuation and for a time-varying state-delay. This was then extended in~\cite{kang2017boundary} for the boundary feedback stabilization of a cascade PDE-ODE system under either Dirichlet or Neumann actuation. In the context of the reaction-diffusion equation, the contribution of this paper is twofold. First, we study the feedback stabilization of the reaction-diffusion equation with Robin-Robin boundary conditions and under a time-varying state-delay. In particular, the proposed approach allows either one single command input (located at one of the two boundaries of the domain) or two command inputs. Second, we show that the resulting closed-loop system is, in the presence of distributed boundary disturbances, exponentially Input-to-State Stable (ISS) with fading memory~\cite{karafyllis2019input}. The concept of ISS, originally introduced by Sontag in~\cite{sontag1989smooth}, plays an important role in the robustness assessment of dynamical systems and the stability of interconnected systems~\cite{karafyllis2019input}. The extension of this notion to infinite-dimensional systems raises many challenges~\cite{mironchenko2016restatements,mironchenko2017characterizations} and is the topic of active research activities. For a complete review of the ISS theory for infinite-dimensional systems, we refer the reader to~\cite{mironchenko2019input}. 

The proposed control design strategy is organized as follows. First, a finite dimensional truncated model is obtained via spectral decomposition. The order of the spectral decomposition is selected to capture all the unstable modes of the original infinite-dimensional system plus a certain number of slow stable modes. In particular, the order is selected to guarantee the robust stability of the residual infinite-dimensional system with respect to exponentially vanishing command inputs exhibiting a prescribed decay rate. We show that this allows the design of the control law based on the finite-dimensional truncated model while ensuring the stability of the full infinite-dimensional closed-loop system.

The remainder of this paper is organized as follows. The problem setting and the control strategy are reported in Section~\ref{sec: problem setting}. The well-posedness of the closed-loop system is assessed in Section~\ref{sec: well-posedness}. The stability analysis is carried out in Section~\ref{sec: stability analysis}. The effectiveness of the proposed control strategy is illustrated by numerical simulations in Section~\ref{sec: numerical simulations}. Concluding remarks are provided in Section~\ref{sec: conclusion}.

\textbf{Notation.} The sets of non-negative integers, positive integers, real, non-negative real, positive real, and complex numbers are denoted by $\mathbb{N}$, $\mathbb{N}^*$, $\mathbb{R}$, $\mathbb{R}_+$, $\mathbb{R}_+^*$, and $\mathbb{C}$, respectively. The real and imaginary parts of a complex number $z$ are denoted by $\operatorname{Re} z$ and $\operatorname{Im} z$, respectively. The field $\mathbb{K}$ denotes either $\mathbb{R}$ or $\mathbb{C}$. The set of $n$-dimensional vectors over $\mathbb{K}$ is denoted by $\mathbb{K}^n$ and is endowed with the Euclidean norm $\Vert x \Vert = \sqrt{x^* x}$. The set of $n \times m$ matrices over $\mathbb{K}$ is denoted by $\mathbb{K}^{n \times m}$ and is endowed with the induced norm denoted by $\Vert\cdot\Vert$. Finally, $\mathrm{AC}_\mathrm{loc}(\mathbb{R}_+;\mathbb{R}^n)$ denotes the set of functions $f : \mathbb{R}_+ \rightarrow \mathbb{R}^n$ that are absolutely continuous on any compact interval of $\mathbb{R}_+$.

\section{Problem setting and control strategy}\label{sec: problem setting}

\subsection{Problem setting}

We are concerned with the boundary feedback stabilization of the following reaction-diffusion equation for Robin boundary conditions and with state-delay:
\begin{subequations}\label{eq: pb setting}
\begin{align}
& y_t(t,x) = a y_{xx}(t,x) + b y(t,x) + c y(t-h(t),x) + d(t,x) \label{eq: pb setting - PDE} \\
& \cos(\theta_1) y(t,0) - \sin(\theta_1) y_x(t,0) = u_1(t) \label{eq: pb setting - BC1} \\
& \cos(\theta_2) y(t,1) + \sin(\theta_2) y_x(t,1) = u_2(t) \label{eq: pb setting - BC2} \\
& y(\tau,x) = \phi(\tau,x), \quad \tau \in [-h_M,0] \label{eq: pb setting - IC} 
\end{align}
\end{subequations}
for $t > 0$ and $x \in (0,1)$. Here we have $a> 0$, $b,c \in \mathbb{R}$ with $c \neq 0$, and $\theta_1,\theta_2 \in [0,2\pi)$. In this setting, $u_1,u_2 : \mathbb{R}_+ \rightarrow \mathbb{R}$ are the boundary controls (with possibly one input set identically equal to zero), $h : \mathbb{R}_+ \rightarrow \mathbb{R}$ is a time-varying delay, and $d : \mathbb{R}_+ \times (0,1) \rightarrow \mathbb{R}$ is a distributed perturbation with the regularity $L^\infty_\mathrm{loc}(\mathbb{R}_+;L^2(0,1))$, i.e. $d : \mathbb{R}_+ \rightarrow L^2(0,1)$ is Bochner measurable and essentially bounded on any compact interval of $\mathbb{R}_+$. We assume that $h$ is continuous and that there exist constants, $0 < h_m < h_M$, such that $h_m \leq h(t) \leq h_M$ for all $t \geq 0$. Finally, $\phi : [-h_M,0] \times (0,1) \rightarrow \mathbb{R}$ represents the initial condition. 

In the sequel, we use the following abstract version of (\ref{eq: pb setting - PDE}-\ref{eq: pb setting - IC}) defined over the state-space $\mathcal{H} = L^2(0,1)$ endowed with its usual inner product $\left< f , g\right> = \int_0^1 f(\xi) g(\xi) \diff\xi$ and associated $L^2$-norm that is also denoted, with a slight abuse of notation, $\Vert \cdot \Vert$.
\begin{subequations}\label{eq: abstract form}
\begin{align}
\dfrac{\mathrm{d}X}{\mathrm{d}t}(t) & = \mathcal{A}X(t) + c X(t-h(t)) + p(t) \label{eq: abstract form - DE} \\
\mathcal{B}X(t) & = u(t) \label{eq: abstract form - BC} \\
X(\tau) & = \Phi(\tau), \quad \tau \in [-h_M,0] \label{eq: abstract form - IC} 
\end{align}
\end{subequations}
for $t>0$ with $\mathcal{A}f = a f'' + b f \in\mathcal{H}$ defined on $D(\mathcal{A}) = H^2(0,1)$, $\mathcal{B}f = (\cos(\theta_1)f(0)-\sin(\theta_1)f'(0),\cos(\theta_2)f(1)+\sin(\theta_2)f'(1)) \in \mathbb{R}^2$ defined on $D(\mathcal{B}) = H^2(0,1)$, $X(t) = y(t,\cdot) \in \mathcal{H}$, $u(t)=(u_1(t),u_2(t)) \in \mathbb{R}^2$, $p(t) = d(t,\cdot) \in \mathcal{H}$, and $\Phi(t)=\phi(t,\cdot) \in \mathcal{H}$. The abstract system (\ref{eq: abstract form}) is a natural extension of the concept of boundary control system reported in~\cite[Def.~3.3.2]{Curtain2012}. Indeed, 1) the disturbance-free operator $\mathcal{A}_0 \triangleq \left.\mathcal{A}\right\vert_{D(\mathcal{A}_0)}$ with $D(\mathcal{A}_0) \triangleq D(\mathcal{A}) \cap \mathrm{ker}(\mathcal{B})$ generates a $C_0$-semigroup denoted by $S(t)$~\cite{delattre2003sturm}; 2) the operator $L_k \in \mathcal{L}(\mathbb{R}^2,\mathcal{H})$, with $k \geq 2$ an integer such that $\cos(\theta_m)+k\sin(\theta_m) \neq 0$ for any $m \in \{1,2\}$, defined for any $u = (u_1,u_2) \in \mathbb{R}^2$ by 
\begin{equation*}
[L_k u](x) = \dfrac{u_1 (1-x)^k}{\cos(\theta_1)+k\sin(\theta_1)} + \dfrac{u_2 x^k}{\cos(\theta_2)+k\sin(\theta_2)}
\end{equation*}
for $x \in [0,1]$, is a lifting operator in the sense that $R(L_k) \subset D(\mathcal{A})$, $\mathcal{A} L_k$ is bounded, and $\mathcal{B}L_k = I_{\mathbb{R}^2}$. Here $\mathrm{ker}(\mathcal{B})$ denotes the kernel of $\mathcal{B}$ while $R(L_k)$ stands for the range of $L_k$.

Assuming that\footnote{This regularity of the forthcoming control law, as well as the existence and uniqueness of the mild solutions for the resulting closed-loop system, will be assessed in Section~\ref{sec: well-posedness}} $u \in \mathrm{AC}_\mathrm{loc}(\mathbb{R}_+;\mathbb{R}^2)$, $p \in L^\infty_\mathrm{loc}(\mathbb{R}_+;\mathcal{H})$, $h \in \mathcal{C}^0(\mathbb{R}_+;\mathbb{R})$ with $0 < h_m \leq h \leq h_M$, and $\Phi \in \mathcal{C}^0([-h_M,0];\mathcal{H})$, the mild solution $X \in \mathcal{C}^0(\mathbb{R}_+;\mathcal{H})$ of (\ref{eq: abstract form}) is uniquely defined by (see~\cite[Def.~3.1.4, Lem.~3.1.5, and Sec.~3.3]{Curtain2012})
\begin{align}
& X(t) 
= S(t) \{ \Phi(0) - L_k u(0) \} + L_k u(t) \label{eq: def mild solution} \\
& \, + \int_0^t S(t-s) \{ \mathcal{A}L_k u(s) - L_k \dot{u}(s) + c X(s-h(s)) + p(s) \} \diff s \nonumber
\end{align}
for $t \geq 0$ and with the initial condition $X(\tau) = \Phi(\tau)$ for all $\tau\in[-h_M,0]$.

\subsection{Proposed control strategy}
Introducing the operator $\mathcal{A}_c \triangleq \mathcal{A} + c I_\mathcal{H}$ defined on $D(\mathcal{A}_c)=D(\mathcal{A})$, (\ref{eq: abstract form - DE}) is equivalent to:
\begin{equation*}
\dfrac{\mathrm{d}X}{\mathrm{d}t}(t) = \mathcal{A}_c X(t) + c \left\{ X(t-h(t)) - X(t) \right\} + p(t).
\end{equation*}
We also introduce $\mathcal{A}_{c,0} = \mathcal{A}_0 + c I_\mathcal{H}$, defined on $D(\mathcal{A}_{c,0}) = D(\mathcal{A}_0)$, which generates a $C_0$-semigroup $T(t)$. Then, using~\cite[Thm.~3.2.1]{Curtain2012}, the mild solution (\ref{eq: def mild solution}) can be equivalently rewritten in function of $T(t)$ as:
\begin{align}
X(t) 
& = T(t) \{ \Phi(0) - L_k u(0) \} + L_k u(t) \label{eq: def mild solution bis} \\
& \phantom{=}\; + \int_0^t T(t-s) \bigg\{ \mathcal{A}_{c} L_k u(s) - L_k \dot{u}(s) \nonumber \\
& \hspace{2cm} + c \{ X(s-h(s)) - X(s) \} + p(s) \bigg\} \diff s \nonumber
\end{align}
for all $t \geq 0$. From the Sturm-Liouville theory (see, e.g., \cite[Sec.~8.6]{renardy2006introduction}), it is well known that $\mathcal{A}_{c,0}$ is a self-adjoint operator whose eigenvalues $(\lambda_n)_{n \geq 1}$ are all real and can be sorted to form a strictly decreasing sequence with $\lambda_n \rightarrow - \infty$ when $n \rightarrow + \infty$. Furthermore, denoting by $e_n$ a unit eigenvector of $\mathcal{A}_{c,0}$ associated with $\lambda_n$, $(e_n)_{n \geq 1}$ forms a Hilbert basis of $\mathcal{H} = L^2(0,1)$.

Introducing the notation $x_n(t) = \left< X(t) , e_n \right>$ we have\footnote{The convergence of the series holds in the norm of the state-space $\mathcal{H} = L^2(0,1)$, i.e. in $L^2$-norm.} $X(t) = \sum\limits_{n \geq 1} x_n(t) e_n$ and $\Vert X(t) \Vert^2 = \sum\limits_{n \geq 1} \vert x_n(t) \vert^2$. The projection of (\ref{eq: def mild solution bis}) onto the Hilbert basis $(e_n)_{n \geq 1}$ with the use of $T(t)z = \sum\limits_{n \geq 1} e^{\lambda_n t} \left< z , e_n \right> e_n$ (see~\cite[Thm.~2.3.5]{Curtain2012} and also~\cite{mironchenko2019local}) and an integration by parts shows that
\begin{align*}
x_{n}(t)
& = e^{\lambda_{n}t} x_{n}(0) + \int_0^t e^{\lambda_{n}(t-s)} g_n(s) \diff s
\end{align*}
with $g_n(t) = \big< - \lambda_{n} L_k u(t) + \mathcal{A}_{c} L_k u(t) + c \{ X(t-h(t)) - X(t) \} + p(t) , e_{n} \big>$. As $g_n$ is integrable on any compact interval, we have $x_{n}\in\mathrm{AC}_\mathrm{loc}(\mathbb{R}_+;\mathbb{R})$ and the following ODE (see also~\cite{lhachemi2018iss}) is satisfied for almost all $t \geq 0$:
\begin{align}
\dot{x}_n(t) 
& = \lambda_n x_n(t) + c \left\{ x_n(t-h(t)) - x_n(t) \right\} \label{eq: prel ODE x_n} \\
& \phantom{=}\; - \lambda_n \left< L_k u(t) , e_n \right> + \left< \mathcal{A}_c L_k u(t) , e_n \right> + \langle p(t) , e_n \rangle . \nonumber 
\end{align} 
We denote by $\{f_1,f_2\}$ the canonical basis of $\mathbb{R}^2$. Introducing for integers $n \geq 1$ and $m \in \{1,2\}$ the quantities $b_{n,m} = - \lambda_n \left< L_k f_m , e_n \right> + \left< \mathcal{A}_c L_k f_m , e_n \right> \in \mathbb{R}$ and $b_n = \begin{bmatrix} b_{n,1} & b_{n,2} \end{bmatrix} \in \mathbb{R}^{1 \times 2}$, we obtain that
\begin{align*}
\dot{x}_n(t) & = \lambda_n x_n(t) + c \left\{ x_n(t-h(t)) - x_n(t) \right\}  \\
& \phantom{=}\; + b_n u(t) + \langle p(t) , e_n \rangle .
\end{align*}
For a given integer $N_0 \geq 1$, which will be discussed in the sequel, the introduction of 
\begin{subequations}
\begin{align}
A & = \mathrm{diag}(\lambda_1,\ldots,\lambda_{N_0}) \in \mathbb{R}^{N_0 \times N_0} , \label{eq: truncated model - A} \\
B & = (b_{n,m})_{1 \leq n \leq N_0 , 1 \leq m \leq 2} \in \mathbb{R}^{N_0 \times 2} , \label{eq: truncated model - B} \\
Y(t) & = \begin{bmatrix} x_1(t) & \ldots & x_{N_0}(t) \end{bmatrix}^\top \in \mathbb{R}^{N_0} , \label{eq: truncated model - Y} \\
D(t) & = \begin{bmatrix} \langle p(t) , e_1 \rangle & \ldots & \langle p(t) , e_{N_0} \rangle \end{bmatrix}^\top \in \mathbb{R}^{N_0} , \label{eq: truncated model - tilde_D} \\
Y_\Phi(\tau) & = \begin{bmatrix} \langle \Phi(\tau) , e_1 \rangle & \ldots & \langle \Phi(\tau) , e_{N_0} \rangle \end{bmatrix}^\top \in \mathbb{R}^{N_0} , \label{eq: truncated model - Y_Phi}
\end{align} 
\end{subequations}
yields $Y\in\mathrm{AC}_\mathrm{loc}(\mathbb{R}_+;\mathbb{R}^{N_0})$ with, for almost all $t \geq 0$,
\begin{subequations}\label{eq: truncated model}
\begin{align}
\dot{Y}(t) & = A Y(t) + c \{ Y(t-h(t)) - Y(t) \} \label{eq: ODE Y} \\
& \phantom{=}\;  + B u(t) + D(t) \nonumber \\
Y(\tau) & = Y_\Phi(\tau), \qquad \tau \in [-h_M,0] \label{eq: ODE Y - IC}
\end{align}
\end{subequations}
The control strategy, inspired by the works~\cite{coron2004global,coron2006global,russell1978controllability} in a delay-free context, relies on the following two steps. First, a feedback control $u = K Y$ is designed to exponentially stabilize the finite-dimensional truncated model (\ref{eq: truncated model}) capturing the unstable dynamics plus an adequate number of slow stable modes of (\ref{eq: pb setting}). This configuration includes the case of a single boundary control input (either at $x=0$ or $x=1$) because one can obtain $u_m = 0$ by setting the $m$-th line of the feedback gain $K$ as $0_{1 \times N_0}$. Specifically, the feedback gain $K \in \mathbb{R}^{2 \times N_0}$ is tuned such that all the poles of $A_\mathrm{cl} = A + BK$ are simple and stable with a sufficiently large decay rate. This procedure is allowed, for either one or two boundary control inputs, by the following Lemma, whose proof is described in Annex~\ref{annex: proof commandability} and where $B_m \in \mathbb{R}^{N_0}$ denotes the $m$-th column of the matrix $B$.

\begin{lem}\label{lem: ODE commandable}
For any given $N_0 \geq 1$, the pairs $(A,B)$, $(A,B_1)$ and $(A,B_2)$ satisfy the Kalman condition~\cite{zhou1996robust}.
\end{lem}

\begin{rem}\label{rem: command indep selection lifing operator}
Note that $b_{n,m}$ is computed based on the selection of a given lifting operator $L_k$. Even if such a lifting operator is not unique, the resulting quantity $b_{n,m}$ is actually independent of the particularly selected lifting operator. See~\cite{lhachemi2019feedback} for details. Consequently, the commandability property stated in Lemma~\ref{lem: ODE commandable} is an intrinsic property of the pair $(\mathcal{A}_c,\mathcal{B})$ associated with (\ref{eq: pb setting}) in the sense that it does not depend on the selection of a particular lifting operator.
\end{rem}

In the second step, we will ensure that the design performed on the finite-dimensional truncated model achieves the exponential stabilization of the full infinite-dimensional system provided the fact that the number of modes $N_0$ used to obtain the truncated model is large enough. This will then lead to the establishment of the following theorem which is stated below.

\begin{thm}\label{thm: main theorem}
Let $0 < h_m < h_M$ be arbitrarily given. Let $N_0 \geq 1$ be such that $\lambda_{N_0 + 1} < - 2 \sqrt{5} \vert c \vert$ and consider the matrices $A$ and $B$ defined by (\ref{eq: truncated model - A}-\ref{eq: truncated model - B}). Let $K \in \mathbb{R}^{2 \times N_0}$ be such that $A_\mathrm{cl} = A + BK$ is Hurwitz with simple eigenvalues $\mu_1 , \ldots , \mu_{N_0} \in \mathbb{C}$ satisfying $\operatorname{Re}\mu_n < - 3 \vert c \vert$ for all $1 \leq n \leq N_0$. Then, there exist constants $\kappa , C_0 , C_1 > 0$ such that, for any initial condition $\Phi \in \mathcal{C}^0([-h_M,0];\mathcal{H})$, any distributed perturbation $p \in L^\infty_\mathrm{loc}(\mathbb{R}_+;\mathcal{H})$, and any delay $h \in \mathcal{C}^0(\mathbb{R}_+;\mathbb{R})$ with $h_m \leq h \leq h_M$, the mild solution $X\in\mathcal{C}^0(\mathbb{R}_+;\mathcal{H})$ of (\ref{eq: abstract form}) with $u = K Y$ satisfies
\begin{align}
\Vert y(t,\cdot) \Vert 
& \leq C_0 e^{- \kappa t} \sup\limits_{\tau \in [-h_M,0]} \Vert \phi(\tau,\cdot) \Vert \label{eq: main theorem - estimate y}\\
& \phantom{\leq}\; + C_1 \esssup\limits_{\tau \in [0,t]} e^{- \kappa (t-\tau)}\Vert d(\tau,\cdot) \Vert \nonumber
\end{align}
for all $t \geq 0$, with control input 
\begin{align}
\Vert u(t) \Vert 
& \leq C_0 \Vert K \Vert e^{- \kappa t} \sup\limits_{\tau \in [-h_M,0]} \Vert \phi(\tau,\cdot) \Vert \label{eq: main theorem - estimate u} \\
& \phantom{\leq}\; + C_1 \Vert K \Vert \esssup\limits_{\tau \in [0,t]} e^{- \kappa (t-\tau)}\Vert d(\tau,\cdot) \Vert , \nonumber
\end{align}
where $y(t,\cdot) = X(t)$, $\phi(t,\cdot) = \Phi(t)$, and $d(t,\cdot) = p(t)$.
\end{thm}

\begin{rem}
The derivation of the design constraints $\operatorname{Re}\mu_n < - 3 \vert c \vert$ and $\lambda_{N_0 + 1} < - 2 \sqrt{5} \vert c \vert$ relies on the derivation of small gain arguments toward the proof of Theorem~\ref{thm: main theorem}. The first (resp. second) constraint is used in the proof of Lemma~\ref{lem: stab truncated model} (resp. Lemma~\ref{lem: infinite-dim part negelected in the design}) to ensure the exponential ISS property of the closed-loop truncated model (resp. the residual infinite-dimensional dynamics).
\end{rem}

\begin{rem}
ISS estimate (\ref{eq: main theorem - estimate y}) is said to have fading memory due to the exponential term in the evaluation of the contribution of the perturbation $d$. This term shows that, as time increases, the contribution of past disturbances on the current magnitude of the state trajectory is vanishing exponentially.
\end{rem}

\begin{rem}
Note that the statement of Theorem~\ref{thm: main theorem} is still valid in the case $c = 0$, i.e. in the absence of state-delay. However, a much simpler proof than the one developed here in the case $c \neq 0$, which in particular allows the case of $A_\mathrm{cl}$ with eigenvalues of arbitrary multiplicity, can be given based on a direct integration and estimation of the ODEs of the spectral reduction.
\end{rem}

The remainder of the paper is devoted to the proof of Theorem~\ref{thm: main theorem} and its numerical illustration.

\section{Well-posedness of the closed-loop system}\label{sec: well-posedness}
In this section, we prove the existence and uniqueness of the mild solutions for (\ref{eq: abstract form}) placed in closed loop with the feedback law $u = K Y$. This ensures the validity of the spectral reduction reported in the previous section.

\begin{lem}\label{lem: well-posedness}
Let $0 < h_m < h_M$, $N_0 \geq 1$, and $K \in \mathbb{R}^{2 \times N_0}$ be arbitrary. For any $\Phi \in \mathcal{C}^0([-h_M,0];\mathcal{H})$, $p \in L^\infty_\mathrm{loc}(\mathbb{R}_+;\mathcal{H})$, and $h \in \mathcal{C}^0(\mathbb{R}_+;\mathbb{R})$ with $h_m \leq h \leq h_M$, there exists a unique mild solution $X\in\mathcal{C}^0(\mathbb{R}_+;\mathcal{H})$ of (\ref{eq: abstract form}) with $u = K Y \in \mathrm{AC}_\mathrm{loc}(\mathbb{R}_+;\mathbb{R}^2)$.
\end{lem}

\textbf{Proof.} 
We show first that $X\in\mathcal{C}^0(\mathbb{R}_+;\mathcal{H})$ such that $Y \in \mathrm{AC}_\mathrm{loc}(\mathbb{R}_+;\mathbb{R}^{N_0})$ is a mild solution of (\ref{eq: abstract form}) with $u = K Y$ if and only if $X\in\mathcal{C}^0(\mathbb{R}_+;\mathcal{H})$ and satisfies, for all $t \geq 0$, equation (\ref{eq: def mild solution}) with $u = K \zeta$,
\begin{align}
\zeta(t) & = e^{(A_\mathrm{cl}-cI) t} Y_\Phi(0) \label{eq: aux abstract form - aux input zeta} \\
& \phantom{=}\; + \int_0^t e^{(A_\mathrm{cl}-cI)(t-\tau)} \{ c Y(\tau-h(\tau)) + D(\tau) \} \diff\tau , \nonumber 
\end{align}
and the initial condition $X(\tau) = \Phi(\tau)$ for all $\tau \in [-h_M,0]$. On one hand, if $X$ is a mild solution of (\ref{eq: abstract form}) with $u = KY \in \mathrm{AC}_\mathrm{loc}(\mathbb{R}_+;\mathbb{R}^2)$, then the developments of Section~\ref{sec: problem setting} show that $Y\in\mathrm{AC}_\mathrm{loc}(\mathbb{R}_+;\mathbb{R}^{N_0})$ satisfies (\ref{eq: truncated model}) and thus $\zeta = Y$. On the other hand, assume that $X\in\mathcal{C}^0(\mathbb{R}_+;\mathcal{H})$ satisfies (\ref{eq: def mild solution}) with $u = K \zeta$ and (\ref{eq: aux abstract form - aux input zeta}). We note from (\ref{eq: aux abstract form - aux input zeta}) that $\zeta \in \mathrm{AC}_\mathrm{loc}(\mathbb{R}_+;\mathbb{R}^{N_0})$ and 
\begin{equation*}
\dot{\zeta}(t) = (A_\mathrm{cl}-cI) \zeta(t) + c Y(t-h(t)) + D(t)
\end{equation*}
for almost all $t \geq 0$. Then $u = K \zeta \in \mathrm{AC}_\mathrm{loc}(\mathbb{R}_+;\mathbb{R}^2)$, showing that (\ref{eq: def mild solution}) indeed makes sense  and, reproducing the developments of Section~\ref{sec: problem setting}, 
\begin{equation*}
\dot{Y}(t) = (A-cI) Y(t) + c Y(t-h(t)) + B K \zeta(t) + D(t)
\end{equation*}
for almost all $t \geq 0$. Consequently we have $\dot{\zeta} - \dot{Y} = (A-cI) (\zeta - Y)$ almost everywhere along with the initial condition $\zeta(0) - Y(0) = Y_\Phi(0) - Y_\Phi(0) = 0$. Thus $\zeta = Y \in \mathrm{AC}_\mathrm{loc}(\mathbb{R}_+;\mathbb{R}^{N_0})$, showing that $X$ is a mild solution of (\ref{eq: abstract form}) with $u = KY$.

To conclude, it remains to show the existence and uniqueness of a function $X\in\mathcal{C}^0(\mathbb{R}_+;\mathcal{H})$ satisfying (\ref{eq: def mild solution}) with $u = K \zeta$ and (\ref{eq: aux abstract form - aux input zeta}). From the regularity assumptions and noting that for any $k \geq 0$, $0 \leq t \leq (k+1)h_m$ implies that $- h_M \leq t - h(t) \leq k h_m$, the existence and uniqueness of such a $X\in\mathcal{C}^0(\mathbb{R}_+;\mathcal{H})$ is immediate by a classical steps procedure and~\cite[Lem.~3.1.5]{Curtain2012}.
\qed

\section{Stability analysis}\label{sec: stability analysis}

This section is devoted to the proof of the main result of this paper: namely, the stability of the closed-loop system.

\subsection{Finite-dimensional truncated model}

We first study the problem of state-feedback stabilization of the finite-dimensional truncated system (\ref{eq: truncated model}).

\begin{lem}\label{lem: stab truncated model}
Let $N_0 \geq 1$ and $0 < h_m < h_M$ be arbitrarily given. Consider the matrices $A$ and $B$ defined by (\ref{eq: truncated model - A}-\ref{eq: truncated model - B}). Let $K \in \mathbb{R}^{2 \times N_0}$ be such that $A_\mathrm{cl} = A + BK$ is Hurwitz with simple eigenvalues $\mu_1 , \ldots , \mu_{N_0} \in \mathbb{C}$ such that $\operatorname{Re}\mu_n < - 3 \vert c \vert$ for all $1 \leq n \leq N_0$. Then, there exist constants $\sigma , C_2, C_3 > 0$ such that, for all $Y_\Phi \in \mathcal{C}^0([-h_M,0];\mathbb{R}^{N_0})$, $h \in \mathcal{C}^0(\mathbb{R}_+;\mathbb{R})$ with $h_m \leq h \leq h_M$, and $D \in L^\infty_\mathrm{loc}(\mathbb{R}_+;\mathbb{R}^{N_0})$, the trajectory $Y \in \mathrm{AC}_\mathrm{loc}(\mathbb{R}_+;\mathbb{R}^{N_0})$ of (\ref{eq: truncated model}) with command input $u = KY$ satisfies 
\begin{align}
\Vert Y(t) \Vert 
& \leq C_2 e^{- \sigma t} \sup\limits_{\tau \in [-h_M,0]} \Vert Y_\Phi(\tau) \Vert \label{eq: exp stab finite-dim part} \\
& \phantom{\leq}\; + C_3 \esssup\limits_{\tau \in [0,t]} e^{- \sigma(t-\tau)} \Vert D(\tau) \Vert \nonumber
\end{align}
for all $t \geq 0$.
\end{lem}

The proof is inspired by the small gain-analysis reported in~\cite[Thm.~2.5]{karafyllis2013delay} for the study of the robustness of predictor feedback with respect to delay mismatches. Nevertheless, 1) we provide here a refinement of the estimates by taking advantage of the properties of the matrix $A_\mathrm{cl}$; 2) we consider the contribution of a disturbance input.

\textbf{Proof.}
We define $\alpha =  - \max\limits_{1 \leq n \leq N_0} \operatorname{Re}\mu_n > 3 \vert c \vert$. Let $\sigma \in (0, \alpha)$ be arbitrarily given. As the eigenvalues of $A_\mathrm{cl}$ are simple, there exists $P \in \mathbb{C}^{N_0 \times N_0}$ such that $P A_\mathrm{cl} P^{-1} = \Lambda \triangleq \mathrm{diag}(\mu_1,\ldots,\mu_{N_0})$. Defining $Z=PY \in \mathrm{AC}_\mathrm{loc}(\mathbb{R}_+;\mathbb{R}^{N_0})$, $Z_\Phi=P Y_\Phi \in \mathcal{C}^0([-h_M,0];\mathbb{R}^{N_0})$, and $\hat{D} = P D \in L^\infty_\mathrm{loc}(\mathbb{R}_+;\mathbb{R}^{N_0})$, we obtain that
\begin{equation}\label{eq: ODE Z}
\dot{Z}(t) = \Lambda Z(t) + c \left\{ Z(t-h(t)) - Z(t) \right\} + \hat{D}(t) ,
\end{equation}
for almost all $t \geq 0$ with the initial condition $Z(\tau) = Z_\Phi(\tau)$ for $\tau \in [-h_M,0]$. Defining $v(t) = Z(t) - Z(t - h(t))$ for all $t \geq 0$, an integration shows that
\begin{align*}
v(t) & = ( e^{\Lambda h(t)} - I ) Z(t-h(t)) \\
& \phantom{=}\; + \int_{t-h(t)}^{t} e^{\Lambda (t-\tau)} \{ - c v(\tau) + \hat{D}(\tau) \} \diff\tau ,
\end{align*}
for all $t \geq h_M$. Noting that $\Vert e^{\Lambda \tau} \Vert = e^{- \alpha \tau}$ and $\Vert e^{\Lambda \tau} - I \Vert \leq 2$ for all $\tau \geq 0$, we obtain that
\begin{align*}
\Vert v(t) \Vert & \leq 2 \Vert Z(t-h(t)) \Vert + \vert c \vert \int_{t-h_M}^{t} e^{-\alpha (t-\tau)} \Vert v(\tau) \Vert \diff\tau \\
& \phantom{\leq} + \int_{t-h_M}^{t} e^{-\alpha (t-\tau)} \Vert \hat{D}(\tau) \Vert \diff\tau ,
\end{align*} 
for all $t \geq h_M$. We evaluate the integral terms as follows ($\psi$ denotes either $v$ or $\hat{D}$):
\begin{align*}
& \int_{t-h_M}^{t} e^{-\alpha (t-\tau)} \Vert \psi(\tau) \Vert \diff\tau \\
& \quad\leq e^{-\sigma t} \int_{t-h_M}^{t} e^{-(\alpha-\sigma) (t-\tau)} e^{\sigma\tau} \Vert \psi(\tau) \Vert \diff\tau \\
& \quad\leq e^{-\sigma t} \dfrac{1 - e^{-(\alpha-\sigma)h_M}}{\alpha - \sigma} \esssup\limits_{\tau \in [t-h_M,t]} e^{\sigma\tau} \Vert \psi(\tau) \Vert ,
\end{align*}
which gives, for all $t \geq h_M$,
\begin{align}
\sup\limits_{\tau \in [h_M,t]} e^{\sigma\tau} \Vert v(\tau) \Vert
& \leq 2 e^{\sigma h_M} \sup\limits_{\tau \in [0,t-h_m]} e^{\sigma\tau} \Vert Z(\tau) \Vert \label{eq: prel estimate v} \\
& \phantom{\leq}\; + \vert c \vert \dfrac{1 - e^{-(\alpha-\sigma)h_M}}{\alpha - \sigma} \sup\limits_{\tau \in [0,t]} e^{\sigma\tau} \Vert v(\tau) \Vert \nonumber \\
& \phantom{\leq}\; + \dfrac{1 - e^{-(\alpha-\sigma)h_M}}{\alpha - \sigma} \esssup\limits_{\tau \in [0,t]} e^{\sigma\tau} \Vert \hat{D}(\tau) \Vert . \nonumber
\end{align}
Now, integrating (\ref{eq: ODE Z}) on $[0,t]$, we obtain for all $t \geq 0$:
\begin{align*}
\Vert Z(t) \Vert 
& \leq e^{-\alpha t} \Vert Z_\Phi(0) \Vert
+ \vert c \vert \int_0^t e^{-\alpha (t-\tau)} \Vert v(\tau) \Vert \diff\tau \\
& \phantom{\leq}\; + \int_0^t e^{-\alpha (t-\tau)} \Vert \hat{D}(\tau) \Vert \diff\tau \\
& \leq e^{-\sigma t} \Vert Z_\Phi(0) \Vert 
+ \dfrac{\vert c \vert e^{-\sigma t}}{\alpha-\sigma} \sup\limits_{\tau \in [0,t]} e^{\sigma\tau} \Vert v(\tau) \Vert \\
& \phantom{\leq}\; + \dfrac{e^{-\sigma t}}{\alpha-\sigma} \esssup\limits_{\tau \in [0,t]} e^{\sigma\tau} \Vert \hat{D}(\tau) \Vert ,
\end{align*}
hence
\begin{align}
\sup\limits_{\tau \in [0,t]} e^{\sigma\tau} \Vert Z(\tau) \Vert
& \leq \Vert Z_\Phi(0) \Vert + \dfrac{\vert c \vert}{\alpha-\sigma} \sup\limits_{\tau \in [0,t]} e^{\sigma\tau} \Vert v(\tau) \Vert \nonumber \\
& \phantom{\leq}\; + \dfrac{1}{\alpha-\sigma} \esssup\limits_{\tau \in [0,t]} e^{\sigma\tau} \Vert \hat{D}(\tau) \Vert . \label{eq: prel estimate Z}
\end{align}
Introducing $\delta \geq 0$ defined by
\begin{equation}\label{eq: finite-dim system - small gain}
\delta = \dfrac{\vert c \vert}{\alpha-\sigma} \left\{ 1 - e^{-(\alpha-\sigma)h_M} + 2 e^{\sigma h_M}  \right\} ,
\end{equation}
we obtain from (\ref{eq: prel estimate v}-\ref{eq: prel estimate Z}) that, for all $t \geq h_M$,
\begin{align*}
\sup\limits_{\tau \in [h_M,t]} e^{\sigma\tau} \Vert v(\tau) \Vert
& \leq 
2 e^{\sigma h_M} \Vert Z_\Phi(0) \Vert + \delta \sup\limits_{\tau \in [0,t]} e^{\sigma\tau} \Vert v(\tau) \Vert \\
& \phantom{\leq}\; + \dfrac{\delta}{\vert c \vert} \esssup\limits_{\tau \in [0,t]} e^{\sigma\tau} \Vert \hat{D}(\tau) \Vert .
\end{align*}
Using now the control design constraint $\alpha > 3 \vert c \vert$, a continuity argument in $\sigma = 0$ shows the existence of a $\sigma \in (0,\alpha)$ such that $\delta < 1$. We fix such a $\sigma \in (0,\alpha)$ for the remaining of the proof. Considering separately the cases where the supremum on $[0,t]$ is achieved either on $[0,h_M]$ or $[h_M,t]$, we obtain for all $t \geq h_M$,
\begin{align}
\sup\limits_{\tau \in [0,t]} e^{\sigma\tau} \Vert v(\tau) \Vert
& \leq 
\dfrac{2 e^{\sigma h_M}}{1-\delta} \Vert Z_\Phi(0) \Vert + \sup\limits_{\tau \in [0,h_M]} e^{\sigma\tau} \Vert v(\tau) \Vert \nonumber \\
& \phantom{\leq}\; + \dfrac{\delta}{\vert c \vert (1-\delta)} \esssup\limits_{\tau \in [0,t]} e^{\sigma\tau} \Vert \hat{D}(\tau) \Vert . \label{eq: estimate sup v - part 1}
\end{align}
We now evaluate the second term on the right-hand side of the above inequality. To do so, we note that, for $0 \leq \tau \leq t \leq h_M$, $\Vert v(\tau) \Vert \leq \Vert Z(\tau) \Vert + \Vert Z(\tau-h(\tau)) \Vert \leq 2 \sup\limits_{\tau \in [0,t]} \Vert Z(\tau) \Vert + \sup\limits_{\tau \in [-h_M,0]} \Vert Z_\Phi(\tau) \Vert$. Furthermore, integrating (\ref{eq: ODE Z}), one can show (using e.g. Gr{\"o}nwall's inequality) the existence of constants $\gamma_0,\gamma_1 > 0$, independent of $Z_\Phi$, $h$, and $\hat{D}$, such that $\Vert Z(t) \Vert \leq \gamma_0 \sup\limits_{\tau \in [-h_M,0]} \Vert Z_\Phi(\tau) \Vert + \gamma_1 \esssup\limits_{\tau \in [0,t]} \Vert \hat{D}(\tau) \Vert$ for all $0 \leq t \leq h_M$. This yields, for all $0 \leq t \leq h_M$,
\begin{align}
\sup\limits_{\tau \in [0,t]} e^{\sigma\tau} \Vert v(\tau) \Vert
& \leq 
(2 \gamma_0 + 1) e^{\sigma h_M} \sup\limits_{\tau \in [-h_M,0]} \Vert Z_\Phi(\tau) \Vert \nonumber \\
& \phantom{\leq}\; + 2 \gamma_1 e^{\sigma h_M} \esssup\limits_{\tau \in [0,t]} e^{\sigma\tau} \Vert \hat{D}(\tau) \Vert . \label{eq: estimate sup v - part 2}
\end{align} 
Combining (\ref{eq: estimate sup v - part 1}-\ref{eq: estimate sup v - part 2}), we obtain that
\begin{align*}
\sup\limits_{\tau \in [0,t]} e^{\sigma\tau} \Vert v(\tau) \Vert
& \leq 
\gamma_2 \sup\limits_{\tau \in [-h_M,0]} \Vert Z_\Phi(\tau) \Vert \\
& \phantom{\leq}\; + \gamma_3 \esssup\limits_{\tau \in [0,t]} e^{\sigma\tau} \Vert \hat{D}(\tau) \Vert .
\end{align*} 
for all $t \geq 0$ with $\gamma_2 = ( 2 \gamma_0 + 1 + 2/(1-\delta) ) e^{\sigma h_M}$ and $\gamma_3 = 2 \gamma_1 e^{\sigma h_M} + \delta / (\vert c \vert (1-\delta))$. Using finally (\ref{eq: prel estimate Z}), and the facts that $\Vert Y \Vert \leq \Vert P^{-1} \Vert \Vert Z \Vert$, $\Vert Z_\Phi \Vert \leq \Vert P \Vert \Vert Y_\Phi \Vert$, and $\Vert \hat{D} \Vert \leq \Vert P \Vert \Vert D \Vert$, the claimed conclusion holds true.
\qed

\begin{rem}
If $K \in \mathbb{R}^{2 \times N_0}$ is selected such that $A_\mathrm{cl} = A + BK$ has real and simple eigenvalues $\mu_1 , \ldots , \mu_{N_0} \in \mathbb{R}$, then the conclusion of Lemma~\ref{lem: stab truncated model} still holds true under the relaxed assumption $\mu_n < - 2 \vert c \vert$ for all $1 \leq n \leq N_0$. This follows from the fact that, in the corresponding proof, the estimate $\Vert e^{\Lambda \tau} - I \Vert \leq 2$ can be in this case replaced by $\Vert e^{\Lambda \tau} - I \Vert \leq 1$ where $\tau \geq 0$.
\end{rem}

\subsection{Infinite-dimensional part of the system neglected in the control design}

We now investigate the robust stability property of the infinite-dimensional part of the system that has been hitherto neglected in our control design. Specifically, the following result holds.

\begin{lem}\label{lem: infinite-dim part negelected in the design}
Let $0 < h_m < h_M$ and $\sigma,C_4,C_5>0$ be arbitrarily given. Let $N_0 \geq 1$ be such that $\lambda_{N_0 + 1} < - 2 \sqrt{5} \vert c \vert$. Then, there exist constants $\kappa\in(0,\sigma)$ and $C_6 , C_7 > 0$ such that, for all $\Phi \in \mathcal{C}^0([-h_M,0];\mathcal{H})$, $p \in L^\infty_\mathrm{loc}(\mathbb{R}_+;\mathcal{H})$, $h \in \mathcal{C}^0(\mathbb{R}_+;\mathbb{R})$ such that $h_m \leq h \leq h_M$, and $u \in \mathrm{AC}_\mathrm{loc}(\mathbb{R}_+;\mathbb{R}^2)$ such that
\begin{align}
\Vert u(t) \Vert + \Vert \dot{u}(t) \Vert 
& \leq C_4 e^{-\sigma t} \sup\limits_{\tau \in [-h_M,0]} \Vert \Phi(\tau) \Vert \label{eq: assumed estimate on u} \\
& \phantom{\leq}\; + C_5 \esssup\limits_{\tau \in [0,t]} e^{- \sigma(t-\tau)} \Vert p(\tau) \Vert \nonumber
\end{align}
for almost all $t \geq 0$, the mild solution $X\in\mathcal{C}^0(\mathbb{R}_+;\mathcal{H})$ of (\ref{eq: abstract form}) satisfies
\begin{align}
\sum\limits_{n \geq N_0 + 1} \vert x_n(t) \vert^2 
& \leq C_6 e^{- 2 \kappa t} \sup\limits_{\tau \in [-h_M,0]} \Vert \Phi(\tau) \Vert^2 \label{eq: exp stab infinite-dim part} \\
& \phantom{\leq}\; + C_7 \esssup\limits_{\tau \in [0,t]} e^{- 2 \kappa (t-\tau)} \Vert p(\tau) \Vert^2 \nonumber 
\end{align}
for all $t \geq 0$, where $x_n(t) = \left< X(t) , e_n \right>$.
\end{lem}

\textbf{Proof.}
Introducing $\beta = - \lambda_{N_0 + 1} / 2 > \sqrt{5} \vert c \vert > 0$, let $\kappa \in (0,\min(\beta,\sigma))$ be arbitrarily given. Then $\lambda_n \leq - 2\beta < -2\kappa < 0$ for all $n \geq N_0 + 1$. We introduce, for $t \geq 0$, $z_n(t) = \langle X(t) - L_k u(t) , e_n \rangle$, $p_n(t) = \langle p(t) , e_n \rangle$, 
\begin{align*}
q_n(t) 
& = \langle \mathcal{A}_c L_k u(t) , e_n \rangle - \langle L_k \dot{u}(t) , e_n \rangle \\
& \phantom{=}\; + c \langle L_k [\chi u](t-h(t)) , e_n \rangle - c \langle L_k u(t) , e_n \rangle ,
\end{align*}
where $\chi$ is defined as the characteristic function of the set $[0,+\infty)$. We also introduce the following series: $Z(t) = \sum\limits_{n \geq N_0 + 1} \vert z_n(t) \vert^2 \leq \Vert X(t) - L_k u(t) \Vert^2$, $P(t) = \sum\limits_{n \geq N_0 + 1} \vert p_n(t) \vert^2 \leq \Vert p(t) \Vert^2$, and $Q(t) = \sum\limits_{n \geq N_0 + 1} \vert q_n(t) \vert^2$ which is such that, for almost all $t \geq 0$,
\begin{align}
Q(t) 
& \leq \gamma_1 e^{- 2 \kappa t} \sup\limits_{\tau\in [-h_M,0]} \Vert \Phi(\tau) \Vert^2 \label{eq: estimate Q} \\
& \phantom{\leq}\; + \gamma_2 \esssup\limits_{\tau\in [0,t]} e^{-2\kappa(t-\tau)} \Vert p(\tau) \Vert^2 , \nonumber
\end{align} 
with $\gamma_1 = 8 C_4^2 \left( \Vert \mathcal{A}_c L_k \Vert^2 + \{ 1 + c^2 ( 1 + e^{2 \kappa h_M} ) \} \Vert L_k \Vert^2 \right)$ and $\gamma_2 = \gamma_1 C_5^2/C_4^2$. 
For $t \geq h_M$, we define $v_n(t) = z_n(t) - z_n(t-h(t))$ and $V(t) = \sum\limits_{n \geq N_0 + 1} \vert v_n(t) \vert^2$. As $z_n \in \mathrm{AC}_\mathrm{loc}(\mathbb{R}_+;\mathbb{R})$, we infer from (\ref{eq: prel ODE x_n}) that, for almost all $t \geq 0$,
\begin{align}
\dot{z}_n(t) 
& = \lambda_n z_n(t) + c \{ [\chi z_n](t-h(t)) - z_n(t) \} \label{eq: EDO z_n} \\
& \phantom{=}\; + c \langle [(1-\chi) \Phi](t-h(t)) , e_n \rangle + q_n(t) + p_n(t) . \nonumber
\end{align}
In particular, we have for almost all $t \geq h_M$
\begin{equation}\label{eq: EDO z_n t geq t_M}
\dot{z}_n(t) = \lambda_n z_n(t) - c v_n(t) + q_n(t) + p_n(t) .
\end{equation} 
In the sequel, we always consider integers $n \geq N_0 + 1$. We introduce for any $t_1 < t_2$ and any real-valued and locally essentially bounded function $\psi$ the notation $\mathcal{I}(\psi,t_1,t_2) = \int_{t_1}^{t_2} e^{- 2 \beta (t_2-\tau)} \vert \psi(\tau) \vert \diff\tau$. Then, the following inequalities hold:
\begin{align*}
\mathcal{I}(\psi,t_1,t_2)
& = e^{-2\kappa t_2} \int_{t_1}^{t_2} e^{-2(\beta-\kappa)(t_2-\tau)} \times e^{2\kappa\tau} \vert \psi(\tau) \vert \diff\tau \\
& \leq e^{-2\kappa t_2} \dfrac{1 - e^{-2(\beta-\kappa)(t_2-t_1)}}{2 (\beta-\kappa)} \esssup\limits_{\tau\in[t_1,t_2]} e^{2\kappa\tau} \vert \psi(\tau) \vert
\end{align*}
and, by Cauchy-Schwarz,
\begin{align*}
\mathcal{I}(\psi,t_1,t_2)^2
& \leq \int_{t_1}^{t_2} e^{-2\beta(t_2-\tau)} \diff\tau \times \mathcal{I}(\psi^2,t_1,t_2) \\
& \leq \dfrac{1 - e^{-2\beta(t_2-t_1)}}{2 \beta} \mathcal{I}(\psi^2,t_1,t_2) .
\end{align*}
For any $t \geq 2 h_M$, we obtain via integration of (\ref{eq: EDO z_n t geq t_M}) over $[t-h(t),t]$ that
\begin{align*}
v_n(t)
& = (e^{\lambda_n h(t)}  - 1) z_n(t-h(t)) \\
& \phantom{=}\; + \int_{t-h(t)}^{t} e^{\lambda_n(t-\tau)} \{ -c v_n(\tau) + q_n(\tau) + p_n(\tau) \} \diff\tau ,
\end{align*}
from which, recalling that $\lambda_n \leq - 2 \beta$, we obtain the estimate
\begin{align*}
\vert v_n(t) \vert^2
& \leq 4 \vert z_n(t-h(t)) \vert^2 + 4 \vert c \vert^2 \mathcal{I}(v_n,t-h(t),t)^2 \\
& \phantom{=}\; + 4 \mathcal{I}(q_n,t-h(t),t)^2 + 4 \mathcal{I}(p_n,t-h(t),t)^2 \\
& \leq 4 \vert z_n(t-h(t)) \vert^2 
+ 4 \vert c \vert^2 \gamma_3 \mathcal{I}(v_n^2,t-h(t),t) \\
& \phantom{\leq}\; + 4 \gamma_3 \mathcal{I}(q_n^2,t-h(t),t) 
+ 4 \gamma_3 \mathcal{I}(p_n^2,t-h(t),t) ,
\end{align*}
with $\gamma_3 = (1-e^{-2\beta h_M})/(2 \beta)$. From (\ref{eq: estimate Q}) we deduce that, for all $t \geq 2 h_M$,
\begin{align*}
V(t) 
& \leq 4 Z(t-h(t)) 
+ 4 \gamma_3 \vert c \vert^2 \mathcal{I}(V,t-h(t),t) \\
& \phantom{\leq}\; + 4 \gamma_3 \mathcal{I}(Q,t-h(t),t)
+ 4 \gamma_3 \mathcal{I}(P,t-h(t),t) \\
& \leq 4 e^{2\kappa h_M} e^{-2\kappa h(t)} Z(t-h(t)) \\
& \phantom{\leq}\; + 4 \gamma_4 \vert c \vert^2 e^{-2\kappa t} \sup\limits_{\tau\in[t-h(t),t]} e^{2\kappa\tau} V(\tau) \\
& \phantom{\leq}\; + 4 \gamma_1 \gamma_4 e^{-2\kappa t} \sup\limits_{\tau\in[-h_M,0]} \Vert \Phi(\tau) \Vert^2 \\
& \phantom{\leq}\; + 4 (1+\gamma_2) \gamma_4 e^{-2\kappa t} \esssup\limits_{\tau\in[0,t]} e^{2 \kappa \tau} \Vert p(\tau) \Vert^2 
\end{align*}
with 
\begin{equation*}
\gamma_4 = \dfrac{(1-e^{-2\beta h_M})(1-e^{-2(\beta-\kappa)h_M})}{4 \beta (\beta-\kappa)} .
\end{equation*}
Consequently, we have for all $t \geq 2 h_M$,
\begin{align}
\sup\limits_{\tau\in[2 h_M,t]} e^{2\kappa\tau} V(\tau)
& \leq 4 e^{2\kappa h_M} \sup\limits_{\tau\in[h_M,t-h_m]} e^{2\kappa\tau} Z(\tau) \nonumber \\
& \phantom{\leq}\; + 4 \gamma_4 \vert c \vert^2 \sup\limits_{\tau\in[h_M,t]} e^{2\kappa\tau} V(\tau) \nonumber \\
& \phantom{\leq}\; + 4 \gamma_1 \gamma_4 \sup\limits_{\tau\in[-h_M,0]} \Vert \Phi(\tau) \Vert^2 \nonumber \\
& \phantom{\leq}\; + 4 (1+\gamma_2) \gamma_4 \esssup\limits_{\tau\in[0,t]} e^{2 \kappa \tau} \Vert p(\tau) \Vert^2 . \label{eq: prel-prel estimate sup V}
\end{align}
Now, integrating (\ref{eq: EDO z_n t geq t_M}) over $[h_M,t]$, we have for $t \geq h_M$
\begin{align*}
\vert z_n(t) \vert
& \leq e^{-2 \beta (t-h_M)} \vert z_n (h_M) \vert + \vert c \vert \mathcal{I}(v_n,h_M,t) \\
& \phantom{\leq}\; + \mathcal{I}(q_n,h_M,t) + \mathcal{I}(p_n,h_M,t) , 
\end{align*}
from which we infer
\begin{align*}
\vert z_n(t) \vert^2
& \leq 4 e^{-4 \kappa (t-h_M)} \vert z_n (h_M) \vert^2 + \dfrac{2 \vert c \vert^2}{\beta} \mathcal{I}(v_n^2,h_M,t) \\
& \phantom{\leq}\; + \dfrac{2}{\beta} \mathcal{I}(q_n^2,h_M,t) + \dfrac{2}{\beta} \mathcal{I}(p_n^2,h_M,t) .
\end{align*}
Thus, using again (\ref{eq: estimate Q}), we obtain for all $t \geq h_M$
\begin{align*}
Z(t)
& \leq 4 e^{- 2 \kappa (t-h_M)} Z(h_M) 
+ \dfrac{2 \vert c \vert^2}{\beta} \mathcal{I}(V,h_M,t) \\
& \phantom{\leq}\; + \dfrac{2}{\beta} \mathcal{I}(Q,h_M,t) + \dfrac{2}{\beta} \mathcal{I}(P,h_M,t) \\
& \leq 4 e^{2 \kappa h_M} e^{- 2 \kappa t} Z(h_M) 
+ \gamma_5 \vert c \vert^2 e^{-2\kappa t} \sup\limits_{\tau\in[h_m,t]} e^{2 \kappa\tau} V(\tau) \\
& \phantom{\leq}\; + \gamma_1 \gamma_5 e^{-2\kappa t} \sup\limits_{\tau\in[-h_M,0]} \Vert \Phi(\tau) \Vert^2 \\
& \phantom{\leq}\; + (1+\gamma_2) \gamma_5 e^{-2\kappa t} \esssup\limits_{\tau\in[0,t]} e^{2 \kappa \tau} \Vert p(\tau) \Vert^2 ,
\end{align*}
with $\gamma_5 = 1/(\beta(\beta-\kappa))$.
Hence, we have for all $t \geq h_M$,
\begin{align}
\sup\limits_{\tau\in[h_M,t]} e^{2\kappa\tau} Z(\tau)
& \leq 4 e^{2 \kappa h_M} Z(h_M) 
+ \gamma_5 \vert c \vert^2 \sup\limits_{\tau\in[h_M,t]} e^{2\kappa\tau} V(\tau) \nonumber \\
& \phantom{\leq}\; + \gamma_1 \gamma_5 \sup\limits_{\tau\in[-h_M,0]} \Vert \Phi(\tau) \Vert^2 \label{eq: prel estimate sup Z} \\
& \phantom{\leq}\; + (1+\gamma_2) \gamma_5 \esssup\limits_{\tau\in[0,t]} e^{2 \kappa \tau} \Vert p(\tau) \Vert^2 . \nonumber
\end{align}
Introducing $\eta \geq 0$ defined by
\begin{align}
\eta & = 4 \vert c \vert^2 ( \gamma_4 + \gamma_5 e^{2 \kappa h_M} ) \label{eq: infinite-dim system - small gain} \\
& = \dfrac{\vert c \vert^2}{\beta (\beta-\kappa)} \left\{ (1-e^{-2\beta h_M})(1-e^{-2(\beta-\kappa)h_M}) + 4 e^{2 \kappa h_M} \right\} , \nonumber
\end{align}
we obtain from (\ref{eq: prel-prel estimate sup V}-\ref{eq: prel estimate sup Z}) that, for all $t \geq 2 h_M$, 
\begin{align*}
\sup\limits_{\tau\in[2 h_M,t]} e^{2\kappa\tau} V(\tau)
& \leq 16 e^{4 \kappa h_M} Z(h_M) + \eta \sup\limits_{\tau\in[h_M,t]} e^{2\kappa\tau} V(\tau) \\
& \phantom{\leq}\; + \dfrac{\gamma_1 \eta}{\vert c \vert^2} \sup\limits_{\tau\in[-h_M,0]} \Vert \Phi(\tau) \Vert^2 \\
& \phantom{\leq}\; + \dfrac{(1+\gamma_2) \eta}{\vert c \vert^2} \esssup\limits_{\tau\in[0,t]} e^{2 \kappa \tau}  \Vert p(\tau) \Vert^2 .
\end{align*}
Using now the control design constraint
$\beta = - \lambda_{N_0+1} / 2 > \sqrt{5} \vert c \vert$, a continuity argument in $\kappa = 0$ shows the existence of a $\kappa \in (0,\min(\beta,\sigma))$ such that $\eta < 1$. We fix such a $\kappa \in (0,\min(\beta,\sigma))$ for the remaining of the proof. As all the supremums in the latter estimate are finite, we obtain for all $t \geq 2 h_M$ that
\begin{align}
\sup\limits_{\tau\in[2 h_M,t]} e^{2\kappa\tau} V(\tau)
& \leq \dfrac{16 e^{4 \kappa h_M}}{1-\eta} Z(h_M) \label{eq: prel estimate sup V} \\
& \phantom{\leq}\; + \dfrac{\eta}{1-\eta} \sup\limits_{\tau\in[h_M,2 h_M]} e^{2\kappa\tau} V(\tau) \nonumber \\
& \phantom{\leq}\; + \dfrac{\gamma_1 \eta}{\vert c \vert^2 (1-\eta)} \sup\limits_{\tau\in[-h_M,0]} \Vert \Phi(\tau) \Vert^2  \nonumber \\
& \phantom{\leq}\; + \dfrac{(1+\gamma_2) \eta}{\vert c \vert^2 (1-\eta)} \esssup\limits_{\tau\in[0,t]} e^{2 \kappa \tau} \Vert p(\tau) \Vert^2 . \nonumber
\end{align}
Based on the integration of (\ref{eq: EDO z_n}), straightforward estimations and the use of (\ref{eq: estimate Q}) show the existence of constants $\gamma_6,\gamma_7 > 0$, independent of $\Phi$, $p$, $h$, and $u$, such that
\begin{equation}\label{eq: rough estimate Z over 0 - h_m}
Z(t) 
\leq \gamma_6 \sup\limits_{\tau\in [-h_M,0]} \Vert \Phi(\tau) \Vert^2    
+ \gamma_7 \esssup\limits_{\tau\in [0,t]} e^{2 \kappa \tau} \Vert p(\tau) \Vert^2, 
\end{equation}
for all $0 \leq t \leq 2 h_M$. Now, noting that, for $h_M \leq \tau \leq t \leq 2 h_M$, $V(\tau) \leq 2 Z(\tau) + 2 Z(\tau-h(\tau))$, we infer that 
\begin{align}
\sup\limits_{\tau\in[h_M,t]} e^{2\kappa\tau} V(\tau) 
& \leq 4 e^{4\kappa h_M} \sup\limits_{\tau\in[0,t]} Z(\tau) \nonumber \\
& \leq 4 \gamma_6 e^{4\kappa h_M} \sup\limits_{\tau\in [-h_M,0]} \Vert \Phi(\tau) \Vert^2 \label{eq: rough estimate V over h_m - 2h_m} \\
& \phantom{\leq}\; + 4 \gamma_7 e^{4\kappa h_M} \esssup\limits_{\tau\in [0,t]} e^{2 \kappa \tau} \Vert p(\tau) \Vert^2 , \nonumber
\end{align}
for all $h_M \leq t \leq 2 h_M$. Combining (\ref{eq: prel estimate sup V}-\ref{eq: rough estimate V over h_m - 2h_m}), we obtain the existence of constants $\gamma_8, \gamma_9 > 0$, independent of $\Phi$, $p$, $h$, and $u$, such that
\begin{align}
\sup\limits_{\tau\in[h_M,t]} e^{2\kappa\tau} V(\tau)
& \leq \gamma_8 \sup\limits_{\tau\in [-h_M,0]} \Vert \Phi(\tau) \Vert^2 \label{eq: estimate V} \\
& \phantom{\leq}\; + \gamma_9 \esssup\limits_{\tau\in [0,t]} e^{2 \kappa \tau}  \Vert p(\tau) \Vert^2 , \nonumber 
\end{align}
for all $t \geq h_M$. Using (\ref{eq: estimate V}) into (\ref{eq: prel estimate sup Z}) and combining with (\ref{eq: rough estimate Z over 0 - h_m}), we obtain the existence of constants $\gamma_{10},\gamma_{11} > 0$, independent of $\Phi$, $p$, $h$, and $u$, such that
\begin{align}
\sup\limits_{\tau\in[0,t]} e^{2\kappa\tau} Z(\tau)
& \leq \gamma_{10} \sup\limits_{\tau\in [-h_M,0]} \Vert \Phi(\tau) \Vert^2  \label{eq: estimate sup Z} \\
& \phantom{\leq}\; + \gamma_{11} \esssup\limits_{\tau\in [0,t]} e^{2 \kappa \tau}  \Vert p(\tau) \Vert^2 , \nonumber
\end{align}
for all $t \geq 0$. Finally, we note that $\vert x_n(t) \vert \leq \vert z_n(t) \vert + \vert \langle L_k u(t) , e_n \rangle \vert$ hence, for all $t \geq 0$,
\begin{align*}
\sum\limits_{n \geq N_0 + 1} \vert x_n(t) \vert^2 
& \leq 2 Z(t) + 2 \Vert L_k \Vert^2 \Vert u(t) \Vert^2 .
\end{align*}
As $u$ is continuous, the right-hand side of (\ref{eq: assumed estimate on u}) is actually an upper-bound of $\Vert u(t) \Vert$ for all $t \geq 0$. Recalling that $0 < \kappa < \sigma$ and using (\ref{eq: estimate sup Z}) the proof is complete.
\qed

\subsection{Proof of the main result and complementary remark}

We can now complete the proof of the main result.

\textbf{Proof of Theorem~\ref{thm: main theorem}.} 
Let $0 < h_m < h_M$ be arbitrarily given. As $(\lambda_n)_{n \geq 1}$ is strictly decreasing with $\lambda_n \rightarrow - \infty$ when $n \rightarrow + \infty$, let $N_0 \geq 1$ be an integer such that $\lambda_{N_0 + 1} < - 2 \sqrt{5} \vert c \vert$. Considering the matrices $A$ and $B$ defined by (\ref{eq: truncated model - A}-\ref{eq: truncated model - B}), let $K \in \mathbb{R}^{2 \times N_0}$ be a feedback gain such that $A+BK$ is Hurwitz with simple eigenvalues $\mu_1 , \ldots , \mu_{N_0} \in \mathbb{C}$ such that $\operatorname{Re}\mu_n < - 3 \vert c \vert$ for all $1 \leq n \leq N_0$. The existence of such a feedback gain is ensured by Lemma~\ref{lem: ODE commandable}. From Lemma~\ref{lem: stab truncated model}, we have the existence of constants $\sigma , C_2, C_3 > 0$ such that, for all $\Phi \in \mathcal{C}^0([-h_M,0]\mathcal{H})$, $p \in L^\infty_\mathrm{loc}(\mathbb{R}_+;\mathcal{H})$, $h \in \mathcal{C}^0(\mathbb{R}_+;\mathbb{R})$ with $h_m \leq h \leq h_M$, the mild solution $X\in\mathcal{C}^0(\mathbb{R}_+\mathcal{H})$ of (\ref{eq: abstract form}) with $u = K Y$ satisfies, for all $t \geq 0$,
\begin{align}
\sqrt{ \sum\limits_{n = 1}^{N_0} \vert x_n(t) \vert^2 }
& \leq C_2 e^{- \sigma t} \sup\limits_{\tau \in [-h_M,0]} \Vert \Phi(\tau) \Vert \label{eq: estimate finite dim part} \\
& \phantom{\leq}\; + C_3 \esssup\limits_{\tau \in [0,t]} e^{- \sigma(t-\tau)} \Vert p(\tau) \Vert , \nonumber
\end{align}
where $x_n(t) = \left< X(t) , e_n \right>$. This estimate follows from (\ref{eq: exp stab finite-dim part}) with $Y$ defined by (\ref{eq: truncated model - Y}) and using the fact that, based on (\ref{eq: truncated model - tilde_D}-\ref{eq: truncated model - Y_Phi}) and recalling that $(e_n)_{n \geq 1}$ is a Hilbert basis, $\Vert D(t) \Vert \leq \Vert p(t) \Vert$ and $\Vert Y_\Phi (\tau) \Vert \leq \Vert \Phi(\tau) \Vert$. From $u = KY$ with $Y \in\mathrm{AC}_\mathrm{loc}(\mathbb{R}_+;\mathbb{R}^{N_0})$ satisfying (\ref{eq: truncated model}), we obtain the existence of constants $C_4 , C_5 > 0$, independent of $\Phi$, $p$, and $h$, such that the estimate (\ref{eq: assumed estimate on u}) holds for almost all $t \geq 0$. Thus, the application of Lemma~\ref{lem: infinite-dim part negelected in the design} yields the existence of constants $\kappa\in(0,\sigma)$ and $C_6 , C_7 > 0$, independent of $\Phi$, $p$, and $h$, such that (\ref{eq: exp stab infinite-dim part}) holds for all $t \geq 0$. Consequently, as $0 < \kappa < \sigma$, we obtain from (\ref{eq: exp stab infinite-dim part}) and (\ref{eq: estimate finite dim part}) that
\begin{align*}
\Vert X(t) \Vert
= \sqrt{ \sum\limits_{n \geq 1} \vert x_n(t) \vert^2 }
& \leq C_0 e^{- \kappa t} \sup\limits_{\tau \in [-h_M,0]} \Vert \Phi(\tau) \Vert \\
& \phantom{\leq}\; + C_1 \esssup\limits_{\tau \in [0,t]} e^{- \kappa(t-\tau)} \Vert p(\tau) \Vert \nonumber
\end{align*}
for all $t \geq 0$, with $C_0 = C_2 + \sqrt{C_6}$ and $C_1 = C_3 + \sqrt{C_7}$ which are constants independent of $\Phi$, $p$, and $h$. As $u = KY$, this completes the proof of Theorem~\ref{thm: main theorem}.
\qed

\begin{rem}
Theorem~\ref{thm: main theorem} merely ensures the existence of a decay rate $\kappa > 0$ but does not provide design mechanisms to further constrain its value. Nevertheless, the above proofs can be used to refine the control design procedure in order to enforce any desired value of the decay rate $\kappa > 0$. Indeed, assume that $0 < h_m < h_M$ and $\kappa > 0$ are arbitrarily given. From (\ref{eq: infinite-dim system - small gain}), we have $\eta \rightarrow 0$ when $\beta \rightarrow + \infty$. Recalling that $ \beta = - \lambda_{N_0 + 1} / 2$ with $\lambda_n \rightarrow - \infty$ when $n \rightarrow + \infty$, we can select an integer $N_0 \geq 1$ such that $\beta > \kappa$ is large enough and thus $\eta < 1$. Now let $\sigma > \kappa$ be given. From (\ref{eq: finite-dim system - small gain}) we have $\delta \rightarrow 0$ when $\alpha \rightarrow + \infty$. Recalling that $\alpha =  - \max\limits_{1 \leq n \leq N_0} \operatorname{Re}\mu_n$, we can select the feedback gain $K \in \mathbb{R}^{2 \times N_0}$ such that $\alpha > \sigma$ is large enough and thus $\delta < 1$. Then, applying the same reasoning as in the proof of Theorem~\ref{thm: main theorem}, we obtain the following result. 
\end{rem}

\begin{thm}
Let $0 < h_m < h_M$ and $\kappa > 0$ be arbitrarily given. There exist $N_0 \geq 1$, $K \in \mathbb{R}^{2 \times N_0}$ and $C_0 , C_1 > 0$ such that, for any initial condition $\Phi \in \mathcal{C}^0([-h_M,0];\mathcal{H})$, any distributed perturbation $p \in L^\infty_\mathrm{loc}(\mathbb{R}_+;\mathcal{H})$, and any delay $h \in \mathcal{C}^0(\mathbb{R}_+;\mathbb{R})$ with $h_m \leq h \leq h_M$, the mild solution $X\in\mathcal{C}^0(\mathbb{R}_+;\mathcal{H})$ of (\ref{eq: abstract form}) with $u = K Y$ satisfies (\ref{eq: main theorem - estimate y}), with control input such that (\ref{eq: main theorem - estimate u}), for all $t \geq 0$.
\end{thm}

\section{Numerical illustration}\label{sec: numerical simulations}
In this section we propose a numerical illustration of the result of Theorem~\ref{thm: main theorem}. Considering the case $h_i \triangleq \cot(\theta_i) > 0$ for $i \in \{1,2\}$, standard computations show that the eigenvalues of $\mathcal{A}_{c,0}$ are given by $\lambda_n = b + c - a r_n^2$ for $n \geq 1$, where $(r_n)_{n \geq 1}$ is the increasing sequence formed by the (strictly) positive solutions $r$ of 
\begin{equation*}
\left( h_1 h_2 - r^2 \right) \sin(r) + ( h_1 + h_2 ) r \cos(r) = 0 .
\end{equation*} 
The corresponding unit eigenvectors are given by $e_n = \phi_n / \Vert \phi_n \Vert$ with $\phi_n(x) = r_n \cos(r_n x) + h_1 \sin(r_n x)$.

\begin{figure}
     \centering
     \includegraphics[width=3.5in]{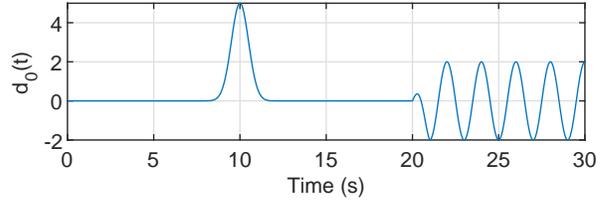}
     \caption{Time evolution of the temporal component $d_0(t)$ of the distributed perturbation $d(t,x)$}
     \label{fig: sim - pert}
\end{figure}

\begin{figure}
     \centering
     	\subfigure[State $y(t,x)$]{
		\includegraphics[width=3.5in]{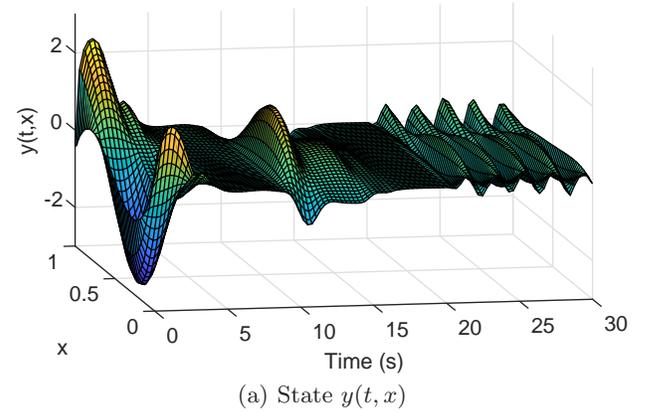}
		\label{fig: sim X - two inputs}
		}
     	\subfigure[Command inputs $u_1$ and $u_2$]{
		\includegraphics[width=3.5in]{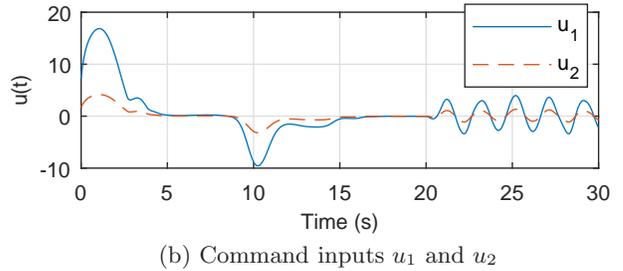}
		\label{fig: sim U - two inputs}
		}
     \caption{Time evolution of the closed-loop system with two command inputs $u_1$ and $u_2$}
     \label{fig: sim - two inputs}
\end{figure}

\begin{figure}
     \centering
     	\subfigure[State $y(t,x)$]{
		\includegraphics[width=3.5in]{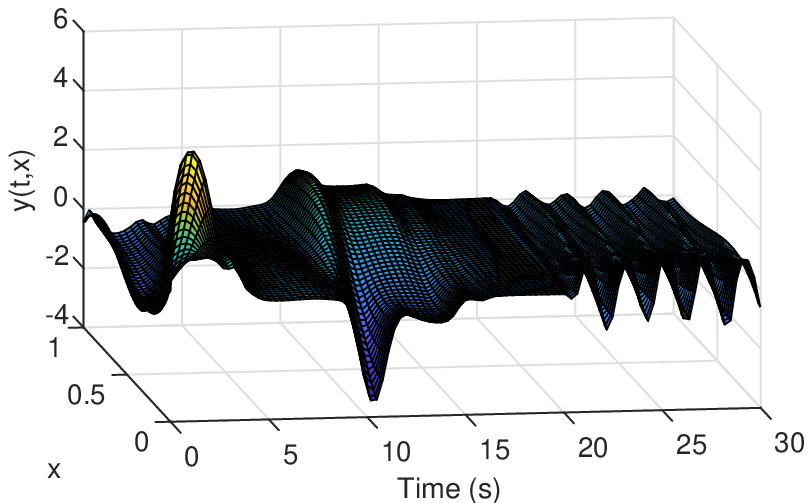}
		\label{fig: sim X - one input}
		}
     	\subfigure[Command input $u_1$]{
		\includegraphics[width=3.5in]{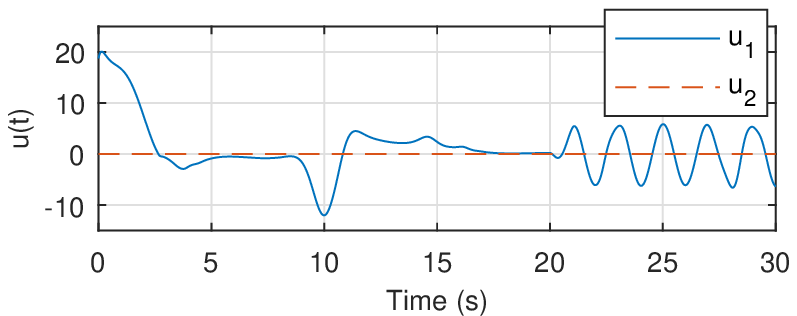}
		\label{fig: sim U - one input}
		}
     \caption{Time evolution of the closed-loop system with one command input $u_1$}
     \label{fig: sim - one input}
\end{figure}

For numerical computations, we set $a = 0.2$, $b = 2$, $c = 1$, $\theta_1 = \pi/3$, and $\theta_2 = \pi/10$. The first three eigenvalues of $\mathcal{A}_{c,0}$ are approximately given by $\lambda_1 \approx 2.5561$, $\lambda_2 \approx -0.1186 > - 2\sqrt{5} \vert c \vert$, and $\lambda_3 \approx -6.2299 < - 2\sqrt{5} \vert c \vert$. Thus we set $N_0 = 2$ and we compute a feedback gain $K \in \mathbb{R}^{2 \times N_0}$ such that the eigenvalues of $A+BK$ are given by $\mu_1 = -3.5$ and $\mu_2 = -4$ with in particular $\mu_2 < \mu_1 < - 3 \vert c \vert$. For simulations, the initial condition, the time-varying delay, and the distributed disturbance are set to $\Phi(t,x) = (1-t)^2 \left\{ (1-2x)/2 + 20 x (1-x) (x-3/5) \right\}$, $h(t) = 2 + 1.5 \sin(t)$, and $d(t,x) = d_0(t) (1-x)$ with $d_0$ as depicted in Fig.\ref{fig: sim - pert}, respectively. The employed numerical scheme relies on the modal approximation of the reaction-diffusion equation using its first 30 modes. The time domain evolution of the closed-loop system (with feedback gain $K$ computed accordingly) is depicted in 1) Fig.~\ref{fig: sim - two inputs} in the case of two command inputs $u_1$ and $u_2$; 2) Fig.~\ref{fig: sim - one input} in the case of a single command input $u_1$ with $u_2 = 0$. As shown over the time interval $[0,8]\,\mathrm{s}$ over which $d(t,\cdot)=0$, both control strategies achieve the exponential stabilization of the closed-loop system with a similar settling time (because the pole placement is identical for the two actuation schemes). Over the time interval $[8,20]\,\mathrm{s}$, the maximum of the perturbation $d$ occurs at time $t = 10\,\mathrm{s}$ and then vanishes as $t$ increases. As expected via the established fading memory estimates (\ref{eq: main theorem - estimate y}-\ref{eq: main theorem - estimate u}), the impact of the vanishing perturbation is rapidly eliminated for $t > 10\,\mathrm{s}$. Finally, for times $t \geq 20\,\mathrm{s}$, we observe the behavior of the closed-loop system in the presence of a non-vanishing disturbance. These results are compliant with the theoretical predictions.

\section{Conclusion}\label{sec: conclusion}
This paper introduced a new method for the feedback stabilization of a reaction-diffusion equation in the presence of a state-delay in the reaction term. The essence of the method relies on the design of the control law on a finite-dimensional truncated model. This LTI model captures the unstable modes and an adequate number of slow stable modes of the original infinite-dimensional system in order to ensure the robust stability of the residual infinite-dimensional dynamics. This technique offers an alternative to backstepping-based methods previously used to tackle this kind of problem. In particular, the developed technique presents the advantage that the feedback is performed on a finite number of modes of the distributed parameter system while the infinite-dimensional residual dynamics is not actively controlled. Future developments of this method include the boundary stabilization of damped wave or beam equations, as well as the derivation of ISS properties with respect to boundary disturbances.


\bibliographystyle{plain}        
\bibliography{autosam}           



\appendix
\section{Proof of Lemma~\ref{lem: ODE commandable}}\label{annex: proof commandability}

As $A$ is diagonal with simple eigenvalues, one can see that the Kalman condition is satisfied provided $b_{n,m} \neq 0$. We show that the latter condition always holds true. Let $k_0 \geq 2$ be such that $\cos(\theta_m) + k \sin(\theta_m) \neq 0$ for all $k \geq k_0$ and all $m \in \{1,2\}$. Following Remark~\ref{rem: command indep selection lifing operator}, we recall that the quantity $b_{n,m}$ is independent of the selected lifting operator $L_k$. Thus, using $\mathcal{A}_c e_n = a e''_n + (b+c)e_n = \lambda_n e_n$ and two successive integration by parts, we have:
\begin{align*}
b_{n,m} & 
= - \lambda_n \left< L_k f_m , e_n \right> + \left< \mathcal{A}_c L_k f_m , e_n \right> \\
& = (b+c - \lambda_n) \left< L_k f_m , e_n \right> + a \left< (L_k f_m)'' , e_n \right> \\
& = a \left\{  - \left< L_k f_m , e''_n \right> + \left< (L_k f_m)'' , e_n \right> \right\} \\
& = a \left\{ (L_k f_m)'(1) e_n(1) - (L_k f_m)(1) e'_n(1) \right\} \\
& \phantom{=} \; + a \left\{ (L_k f_m)(0) e'_n(0) - (L_k f_m)'(0) e_n(0) \right\} ,
\end{align*}
for all $k \geq k_0$, and then 
\begin{align*}
b_{n,1} = a \dfrac{e'_n(0) + k e_n(0)}{\cos(\theta_1) + k \sin(\theta_1)} , \; 
b_{n,2} = a \dfrac{- e'_n(1) + k e_n(1)}{\cos(\theta_2) + k \sin(\theta_2)} .
\end{align*}
Now, the condition $b_{n,1} = 0$ implies $e'_n(0) + k e_n(0) = 0$ for all $k \geq k_0$. This yields $e_n(0) = e'_n(0) = 0$ and we infer by Cauchy uniqueness the contradiction $e_n = 0$. Thus $b_{n,1} \neq 0$ for all $n \geq 1$. The same reasoning shows that $b_{n,2} \neq 0$ for all $n \geq 1$. This concludes the proof.


\end{document}